\documentclass[a4paper,12pt]{article}
\makeatletter

\newcommand{\Rmnum}[1]{\expandafter\@slowromancap\romannumeral #1@}
\makeatother

\makeatletter 
\@addtoreset{equation}{section}
\makeatother  

\makeatletter
\renewcommand{\maketitle}{\bgroup\setlength{\parindent}{0pt}
\begin{flushleft}
  \textbf{\@title}

\vspace{10mm}

  \@author
\end{flushleft}\egroup
}

\usepackage[all]{xy}
\usepackage[depth=subsection]{bookmark}
\usepackage{amsmath}
 \usepackage{relsize}
\usepackage{enumerate}
\usepackage{amsthm}
\usepackage{amsfonts}
\usepackage{dsfont}
\usepackage{bm}
\usepackage{xcolor}
\usepackage{CJK}
\usepackage{txfonts}
\usepackage{multirow}
\usepackage{graphicx}
\usepackage{tikz-cd}
\usepackage{mathtools}

\newtheorem{define}{Definition}[section]
\newtheorem{theorem}[define]{Theorem}
\newtheorem{rem}[define]{Remark}

\newtheorem{cor}[define]{Corollary}
\newtheorem{lem}[define]{Lemma}

\newcommand{\pf}{{\noindent\textbf{Proof}\quad}}

\newcommand{\KAI}{\CJKfamily{kai}}

\def\re{\operatorname{Re}}
\def\spin{\operatorname{spin}}
\def\Aut{\operatorname{Aut}}

\def\mo{\operatorname{mod}}

\def\det{\operatorname{det}}

\def\Tr{\operatorname{Tr}}

\def\SLNZ{\operatorname{SL}(2,\mathbb{Z})}

\def\SpZ{\operatorname{Sp}_{4}(\mathbb{Z})}
\def\triv{\operatorname{triv}}

\title{\LARGE\textbf{Linear Relations of Siegel Poincar\'e Series\\
\vspace{1mm}
and Non-vanshing of the Central Values of \\
\vspace{2mm}
Spinor L-functions \footnote{This version might be different from the version to be published.}}}

\author{Zhining Wei \\
Department of Mathematics, Ohio State University, Columbus, OH 43210, USA\\
wei.863@buckeyemail.osu.edu}

\date{}

\providecommand{\keywords}[1]
{
  \small	
  \noindent\textbf{Keywords} #1
}

\providecommand{\MSC}[1]
{
  \small	
  \noindent\textbf{Mathematics Subject Classification} #1
}

\begin{document}

\maketitle

\vspace{10mm}

\noindent\textbf{Abstract}\hspace{4mm} 
{In this paper, we will first investigate the linear relations of a one parameter family of Siegel Poincar\'e series. Then we give the applications to the non-vanishing of Fourier coefficients of Siegel cusp eigenforms and the central values.
}

\vspace{4mm}

\keywords\hspace{4mm}{Siegel Poincar\'e series, Fourier coefficients of Siegel cusp forms, non-vanishing of central values, B\"ocherer conjecture.}

\vspace{4mm}

\MSC\hspace{4mm}{	11F46, 11F30, 11F67}

\vspace{20mm}

\section{Introduction}
Fourier coefficients of cusp forms are fundamental objects in number theory and there are many open questions. One basic problem is to determine whether a Fourier coefficient vanishes or not. Such problems have a significant application in the theory of L-functions. A remarkable example are the $\operatorname{GL}(2)$ cusp forms. Indeed, due to Waldspurger's theorem, Fourier coefficients of a half-integral weight cusp form will determine the central values of L-functions of the corresponding integral weight cusp form. 

To study the Fourier coefficients of cusp forms, an effective tool are the Poincar\'e series. It is known that the space of holomorphic $\SLNZ$ cusp forms is a vector space of finite dimension. Additionally, the Poincar\'e series can determine Fourier coefficients of cusp forms via the Petersson inner product and hence the space can be generated by Poincar\'e series. Therefore, non-vanishing results and linear relations of Poincar\'e series will reveal the information for Fourier coefficients. 

Such problems have been investigated by several authors. For the non-vanishing of Poincar\'e series, one can refer to \cite{Ra}, \cite{Mo}. For the linear relations of Poincar\'e series,  Petersson's paper \cite{Pe} showed that, the first $m$ Poincar\'e series span the (finite-dimensional) vector space of weight $k$ full level holomorphic cusp forms, where $m$ is the dimension of the space. Later, Petersson's result was generalized to Hecke congruence groups by Lehner \cite{Le} under some conditions.

In this paper, we will consider the situation of Siegel cusp forms. Similar to the $\operatorname{GL}(2)$ case, non-vanishing results and linear relations of Siegel Poincar\'e series are effective tools for studying the Fourier coefficients of Siegel cusp forms, especially Hecke Siegel cusp forms.

The non-vanishing of Siegel Poincar\'e series was studied by several papers. (See \cite{DS}, \cite{DKS}.) In our paper, we will investigate the linear relations of a one parameter family of Siegel Poincar\'e series. 

We proceed to describe our results. Let $\Gamma_2=\SpZ.$  Assume that $k$ is an even number and $F$ a Siegel cusp form of weight $k.$ This is a function on Siegel's half plane $\mathbb{H}_2=\{Z=X+iY\in\operatorname{Mat}_2(\mathbb{C})|Z={}^tZ,\hspace{2mm}Y>0\}$ satisfying
\[F(\gamma Z)=F((AZ+B)(CZ+D)^{-1})=(\det J(\gamma,Z))^{k}F(Z)\hspace{4mm}\mbox{for $\gamma=\begin{pmatrix}A&B\\C&D\end{pmatrix}\in\SpZ$}\]
where $J(\gamma, Z)=CZ+D$ and some holomorphy conditions on $\mathbb{H}_2.$ For a general introduction to Siegel cusp forms, see \cite{Kl}, \cite{Pi}. Denote by $\mathcal{S}_k(\Gamma_2)$ the vector space generated by Siegel cusp forms of weight $k.$ It is known that $\mathcal{S}_k(\Gamma_2)$ has finite dimension. By the automorphy and the holomorphy of $F,$ it has Fourier expansion
\[F(Z)=\sum_{T}A(F,T)e(\Tr(TZ)),\]
where the summation is over all symmetric, positive definite and half-integral matrices $T$. We equip the space $\mathcal{S}_k(\Gamma_2)$ with the Petersson inner product:
\[<F,G>=\int_{\Gamma_2\setminus\mathbb{H}_2}F(Z)\overline{G(Z)}(\det Y)^k\frac{dXdY}{(\det Y)^3}.\]
By the Hecke theory, we can find an (orthogonal) basis, denoted by $\mathcal{B}_k,$ of $\mathcal{S}_k(\Gamma_2)$ which are eigenfunctions of Hecke operators. Then for $F\in\mathcal{B}_k,$ we can define the spinor L-function, denoted by $L(s,\spin, F)$. By some normalization, we can assume that the critical strip is $0<\re s<1.$ (Indeed, this is the spinor L-function of $\pi_F,$ which is the automorphic cuspidal representation associated to $F.$ See Chapter 6 in \cite{Pi}.) 

In addition, the space $\mathcal{S}_k(\Gamma_2)$ contains a subspace of (Saito-Kurokawa) lifts from holomorphic $\SLNZ$ cusp forms of weight $2k-2.$ Denote by $\mathcal{K}$ the subspace of Saito-Kurokawa lifts. Furthermore, let $f$ be an $\SLNZ$ Hecke eigenform of weight $2k-2.$ Then the corresponding Siegel cusp form $F_f$ is also a Hecke eigenform. The spinor L-function of $F_f$ and the (normalized) L-function of $f$ are connected by the following relation:
\[L(s,\spin,F_f)=\zeta(s+1/2)\zeta(s-1/2)L(s,f).\] 

Denote by $\Gamma_{\infty}$ the group of translations in $\Gamma_2,$ that is, 
\[\Gamma_{\infty}=\left\{\left.\begin{pmatrix}I_2&X\\&I_2\end{pmatrix}\in\SpZ\right|X={}^tX\right\}.\]
In the theory of Siegel cusp forms, the Poincar\'e series are defined by
\[\mathbb{P}_{Q}(Z)=\sum_{\gamma\in\Gamma_{\infty}\setminus\SpZ}\det(J(\gamma,Z))^{-k}e(\Tr(Q\gamma Z))=\sum_{T}A(\mathbb{P}_Q,T)e(\Tr(TZ))\]
for any $Q$, a symmetric, positive definite and half-integral matrix. (The explicit form of $A(\mathbb{P}_Q,T)$ will be given in Lemma \ref{3.5}) The following formula explains why Siegel Poincar\'e series can determine Fourier coefficients of Siegel cusp forms (\cite{KST}, (3.1.1)):
\begin{equation}
<F,\mathbb{P}_Q>=8c_k\frac{A(F,Q)}{(\det Q)^{l}}\label{1.1}
\end{equation}
where 
\[l=k-3/2\hspace{10mm}\mbox{and}\hspace{10mm}c_k=\frac{\pi^{1/2}}{4}(4\pi)^{3-2k}\Gamma(k-3/2)\Gamma(k-2).\]

In this paper, we will study Siegel Poincare series $\mathbb{P}_Q$ with $Q=nI_2$ for $n$ a positive integer. We will establish the following theorem:
\begin{theorem}
\upshape\KAI{Let $\epsilon>0$. Denote by $\mathcal{V}$ the vector space spanned by $\mathbb{P}_{nI_2}.$  Then for sufficiently large $k$, we can find a subspace $\mathcal{W}\subseteq\mathcal{V},$ such that
\[\dim\mathcal{W}\geq k^{2/3-\epsilon},\]
and
\[\mathcal{W}\cap\mathcal{K}=\{0\}.\] 
}\label{2.1}
\end{theorem}
In section 3, the construction of the subspace $\mathcal{W}$ will be given explicitly. To prove Theorem \ref{2.1}, we need to investigate the matrix $\mathcal{A}$ defined in \ref{4.1}. We will show that the matrix $\mathcal{A}$ is invertible and this will give the linear independence of Siegel Poincar\'e series. For the second part, we need the Maass relation \ref{5.1}.

As an application, we will study the number of Hecke eigenforms $F$ satisfying $A(F,I_2)\neq 0$ in Section 4. We can establish the following result:  
\begin{cor}
\upshape\KAI{Let $\epsilon>0$ be arbitrary. Then for sufficiently large $k,$
\[\#\{F\in\mathcal{B}_k|\mbox{$F$ is not a Saito-Kurokawa lift and $A(F,I_2)\neq0$}\}\geq k^{2/3-\epsilon}.\]
}\label{2.2}
\end{cor}

In \cite{Bo}, B\"ocherer made a remarkable conjecture that $A(F,I_2)$ should be related to central L-values. This can be regarded as a generalization of Waldspurger's theorem. A precise version of B\"ocherer's conjecture was given in \cite{DPSS}. In our case, the statement is, for a non-Saito-Kurokawa lift $F$,
\[\frac{\pi^{1/2}}{4}(4\pi)^{3-2k}\Gamma(k-3/2)\Gamma(k-2)\frac{A(F,I_2)^2}{||F||^2}=\frac{64\pi^6\Gamma(2k-4)}{\Gamma(2k-1)}\frac{L(1/2,\spin,F)L(1/2,\spin,F\times\chi_{-4})}{L(1,\pi_F,\operatorname{Ad})}\]
where $L(s,\pi_F,\operatorname{Ad})$ is a degree $10$ L-function. By \cite{PSS} Theorem 5.2.1, we know that $L(1,\pi_F,\operatorname{Ad})\neq 0.$ Later B\"ocherer's conjecture was proved by Furusawa and Morimoto \cite{FM}. In our case, we have the following corollary:
\begin{cor}
\upshape\KAI{
Let $\epsilon>0$ be arbitrary. Then for sufficiently large $k,$
\[\#\{F\in\mathcal{B}_k|\mbox{$F$ is not a Saito-Kurokawa lift and $L(1/2,\spin,F)L(1/2,\spin,F\times\chi_{-4})\neq 0$}\}\geq k^{2/3-\epsilon}.\]
}\label{2.4}
\end{cor}
\begin{rem}
In Blomer's paper \cite{Bl}, he proved that, for sufficiently large $k,$ there exists one Hecke eigenform $F$ (not a Saito-Kurokawa lift) such that 
\[L(1/2,\spin,F)L(1/2,\spin,F\times\chi_{-4})\neq 0.\]
(Corollary 5 with $q_1=q_2=1.$) In Corollary 4 of \cite{Bl}, Blomer established a much stronger quantitative lower bound for the size of the set $(k^{3-\epsilon})$ under the ``no-Siegel-zero'' hypothesis for the adjoint $L$-function. 
\end{rem}

\section{Preliminaries}

In this paper, we will always assume $k$ to be an even integer and $k\geq 6.$ Then we set $l=k-\frac{3}{2}$. We also set
\[c_k=\frac{\pi^{1/2}}{4}(4\pi)^{3-2k}\Gamma(k-3/2)\Gamma(k-2).\]

This paper relies heavily on the estimation of the Fourier coefficients of Siegel Poincar\'e series. We will use some results in Section 3 of Blomer's paper \cite{Bl}.

The explicit form of the Fourier coefficient of Siegel Poincar\'e series was calculated by Kitaoka \cite{Ki}. Here we introduce some notations and quote his results. Let  $Q,T$ be positive definite, symmetric and half-integral matrices. Then for an invertible matrix $C\in\operatorname{Mat}_2(\mathbb{C}),$ we define the symplectic Kloosterman sum
\[K(Q,T;C)=\sum e(\Tr(AC^{-1}Q+C^{-1}DT))\]
where the summation is over matrices$\left\{\left.\begin{pmatrix}A&*\\ C&D\end{pmatrix}\right|\begin{pmatrix}A&*\\ C&D\end{pmatrix}\in\SpZ\right\}$ for a given value $C$ in a system $X(C)$ of representatives for $\Gamma_{\infty}\setminus\SpZ/ \Gamma_{\infty}.$ One can see that $|X(C)|\leq |\det C|^{3/2}$ and hence the sum is trivially bounded by $|\det C|^{3/2}.$ For the details, one can see the argument below (3.1) in  Blomer's paper \cite{Bl}.

Let $P$ be a real and positive definite $2\times 2$ matrix with eigenvalues $s_1^2$ and $s_2^2$, then we define
\[\mathcal{J}_l(P)=\int_0^{\pi/2}J_l(4\pi s_1\sin\theta)J_l(4\pi s_2\sin\theta)\sin\theta\,d\theta.\]
For $P=\begin{pmatrix}p_1&p_2/2\\p_2/2&p_4\end{pmatrix}$ and $S=\begin{pmatrix}s_1&s_2/2\\s_2/2&s_4\end{pmatrix}$ being positive definite, symmetric and half-integral matrices and $c\in\mathbb{N}$, we define another Kloosterman sum
\[H^{\pm}(P,S;c)=\delta_{s_4=p_4}{\sum_{d_1 (\mo c)}}^*\sum_{d_2(\mo c)}e\left(\frac{\overline{d_1}s_4d_2^2\mp\overline{d_1}p_2d_2+s_2d_2+\overline{d_1}p_1+d_1s_1}{c}\mp\frac{p_2s_2}{2cs_4}\right).\]
Here we will only use the trivial bound $|H^{\pm}(P,S;c)|\leq c^2.$

The following lemma is due to Kitaoka \cite{Ki}:
\begin{lem}
\upshape\KAI{Let $k\geq 6$ be even. Then for $P,Q$ being positive definite, symmetric and half-integral matrices,  
\begin{flalign*}
A(\mathbb{P}_Q,T)&=\delta_{Q\sim T}\#\operatorname{Aut}(T)\\
&\hspace{2mm}+\left(\frac{\det T}{\det Q}\right)^{l/2}\sum_{\pm}\sum_{U,V}\sum_{s,c\geq 1}\frac{(-1)^k\sqrt{2}\pi}{c^{3/2}s^{1/2}}H^{\pm}(UQ {}^{t}U,V^{-1}T {}^{t}V^{-1},c)J_{l}\left(\frac{4\pi\sqrt{\det TQ}}{cs}\right)\\
&\hspace{2mm}+8\pi^2\left(\frac{\det T}{\det Q}\right)^{l/2}\sum_{\det C\neq 0}\frac{K(Q,T;C)}{|\det C|^{3/2}}\mathcal{J}_l(TC^{-1}Q {}^{t}C^{-1}).
\end{flalign*}
where the sum over $U,V\in\operatorname{GL}_2(\mathbb{Z})$ in the second term on the right hand side is over matrices
\[U=\begin{pmatrix}*&*\\ u_3&u_4\end{pmatrix}/ \{\pm I_2\}\hspace{4mm}V=\begin{pmatrix}v_1&*\\ v_3&*\end{pmatrix}\hspace{4mm}\begin{pmatrix}u_3&u_4\end{pmatrix}Q\begin{pmatrix}u_3\\ u_4\end{pmatrix}=\begin{pmatrix}-v_3&v_1\end{pmatrix}T\begin{pmatrix}-v_3\\ v_1\end{pmatrix}=s.\]
$Q\sim T$ means the equivalence in the sense of quadratic forms and $\operatorname{Aut}(T)=\{U\in\operatorname{GL}_2(\mathbb{Z})|{}^{t}UTU=T\}.$
}\label{3.5}
\end{lem}
Generally, these three terms are called rank $0,$ rank $1$ and rank $2$ terms respectively.

We also recall a basic lemma in Blomer's paper (\cite{Bl} Section 3, Lemma 2):
\begin{lem}
\upshape\KAI{
For positive definite matrices $T,Q$ with largest eigenvalues $\lambda_T, \lambda_Q$, the smallest eigenvalue of $TC^{-1}Q{}^tC^{-1}$ is $<<\lambda_T\lambda_Q||C||^{-2}$
where $||\cdot||$ denotes any submultiplicative matrix norm.
}\label{3.4}
\end{lem}

Additionally, we need the following uniform bounds for Bessel functions:
\begin{equation}
J_l(x)<<\left(\frac{cx}{l}\right)^l,\label{6.5}
\end{equation}
where $c$ can be chosen to be any constant larger than $\frac{e}{2}.$This is valid for $x>0$ and $l>1/2$. This follows from [\cite{GR}, 8.411.4].

Let $Y=k^{2/3-\epsilon/2}$ and we set $\mathcal{P}$ to be the set of all primes in the interval $[Y,2Y].$
\begin{lem}
\upshape\KAI{
For $p,q\in\mathcal{P},$ we have:
\[A(\mathbb{P}_{pI_2},qI_2)=8\delta_{p,q}+O\left(\frac{M^{l}}{l^{(1/3+\epsilon/2)l}}\right).\]
Here $M$ is a constant independent from $k.$ 
}\label{3.1}
\end{lem}
\pf By Kitaoka's formula, we have
\begin{flalign*}
A(\mathbb{P}_{pI_2},qI_2)&=\delta_{pI_2\sim qI_2}\#\operatorname{Aut}(qI_2)\\
&\hspace{2mm}+\left(\frac{\det (qI_2)}{\det (pI_2)}\right)^{l/2}\sum_{\pm}\sum_{U,V}\sum_{s,c\geq 1}\frac{(-1)^k\sqrt{2}\pi}{c^{3/2}s^{1/2}}H^{\pm}(pU {}^{t}U,qV^{-1} {}^{t}V^{-1},c)J_{l}\left(\frac{4\pi\sqrt{\det (pqI_2)}}{cs}\right)\\
&\hspace{2mm}+8\pi^2\left(\frac{\det (qI_2)}{\det (pI_2)}\right)^{l/2}\sum_{\det C\neq 0}\frac{K(pI_2,qI_2;C)}{|\det C|^{3/2}}\mathcal{J}_l(pqC^{-1}{}^{t}C^{-1}).
\end{flalign*}
Notice that $pI_2\sim qI_2$ if and only if $p=q.$ On the other hand, $\#\Aut(pI_2)=\#\Aut(I_2)=8.$ This handled the rank $0$ term. 

Next, by the definition of the exponential sum in the rank $1$ term, we know that $s=p(u_3^2+u_4^2)=q(v_1^2+v_3^2)$ and hence $[p,q]|s.$ Here $[p,q]$ is the least common multiple. Then by the change of variable $s\mapsto [p,q]s,$ the estimation (\ref{6.5}) and the fact that $p,q\in[Y,2Y]$ we know that the rank $1$ term is bounded by
\[\frac{M^l}{l^{(1/3+\epsilon/2)l}}.\]

Finally, we consider the rank $2$ term. We know that $K(pI_2,qI_2;C)$ is bounded by $|\det(C)|^{3/2}.$ By Lemma \ref{3.4}, the smallest eigenvalue of $pqC^{-1} {}^{t}C^{-1}$ is $<<\frac{pq}{||C||^2},$ where $||C||=\Tr(C{}^tC)$ is a matrix norm on $2\times 2$ matrices. This shows that
\[\mathcal{J}(pqC^{-1} {}^{t}C^{-1})<<\left(\frac{M'(pq)^{1/2}}{||C|| l}\right)^l\]
for some constant $M'$ independent from $k.$ Then the rank $2$ term is bounded by
\[\frac{M^l}{l^{(1/3+\epsilon/2)l}}.\]
\qed

\begin{lem}
\upshape\KAI{
For $p,q\in\mathcal{P},$ we have:
\[A\left(\mathbb{P}_{pI_2},\begin{pmatrix}q^2&\\&1\end{pmatrix}\right)<<\frac{M^{l}}{l^{3\epsilon l/4}}.\]
Here $M$ is a constant independent from $k.$  
}\label{3.2}
\end{lem}
\pf The proof can be quite similar to Lemma \ref{3.1}. Here I would only mention two differences: for the rank $1$ term, we can only do the substitution $s\mapsto ps.$ For the rank $2$ term, the smaller eigenvalue of $p\begin{pmatrix}q^2&\\&1\end{pmatrix}C^{-1}{^tC^{-1}}$ is bounded by 
\[\frac{q^2p}{||C||^2}\]
by Lemma \ref{3.4}.
\qed

We also need the following lemma:
\begin{lem}
\upshape\KAI{
For $p\in\mathcal{P},$ we have:
\[A(\mathbb{P}_{pI_2},I_2)<<\frac{M^l}{l^{(4/3-\epsilon/4)l}}.\]
Here $M$ is a constant independent from $k.$  
}\label{3.3}
\end{lem}
\pf The proof is similar to Lemma \ref{3.1}.
\qed

\section{Proof of Theorem \ref{2.1}}
\pf Let $Y=k^{2/3-\epsilon/2}$ and we set $\mathcal{P}$ to be the set of all primes in the interval $[Y,2Y].$ Set $n=\#\mathcal{P}.$ Then by prime number theory, $n\geq k^{2/3-\epsilon}$ when $k$ is sufficiently large.

We need to show two things: the set $\{\mathbb{P}_{pI_2}|p\in\mathcal{P}\}$ is linearly independent. Then let $\mathcal{W}$ be the subspace of $\mathcal{S}_k(\Gamma_2)$ spanned by $\{\mathbb{P}_{pI_2}|p\in\mathcal{P}\}$ and we need to show that $\mathcal{W}\cap\mathcal{K}=\{0\}.$

\subsection{The linear independence}

To show the linear independence of $\{\mathbb{P}_{pI_2}|p\in\mathcal{P}\},$ we consider the following matrix 
\begin{equation}
\mathcal{A}=\left(a_{p,q}\right)_{p,q\in\mathcal{P}}=\left(A(\mathbb{P}_{pI_2},qI_2)\right)_{p,q\in\mathcal{P}},\label{4.1}
\end{equation} 
that is, $a_{p,q}=A(\mathbb{P}_{pI_2},qI_2),$ the $qI_2$-th Fourier coefficient of $\mathbb{P}_{pI_2}.$ We want to show that the matrix $\mathcal{A}$ is invertible and hence Siegel Poincare series $\{\mathbb{P}_{pI_2}|p\in\mathcal{P}\}$ are linearly independent in $\mathcal{V}.$ This can be obtained by using Lemma \ref{3.1}. Indeed, we can write
\[\mathcal{A}=8I_n+\mathcal{B}\]
such that every entry of $\mathcal{B}$ is uniformly bounded by $M^l/l^{(1/3+\epsilon/2)l}<<\frac{1}{k^{100}}$ when $k$ is sufficiently large. Then $\mathcal{A}$ is invertible.

\begin{rem}
Indeed, by a more delicate estimation, we can show, there exists a subset $\mathcal{Q}$ of natural numbers with $\#\mathcal{Q}>>k$ such that
\[\{\mathbb{P}_{nI_2}|n\in\mathcal{Q}\}\]
is a linearly independent set. We will discuss this in the last section.
\end{rem}

\subsection{The Non-intersection of $\mathcal{W}$ and $\mathcal{K}$}
We prove this by contradiction. Suppose not, we can find $(\lambda_p)_{p\in\mathcal{P}}$ such that 
\[F=\sum_{p\in\mathcal{P}}\lambda_p\mathbb{P}_{pI_2}\in\mathcal{K}\]
and not all $\lambda_p$ are $0.$ Then without loss of generality, we can assume that $|\lambda_p|\leq 1$ for all $p\in\mathcal{P}$ and there exists a $p_0$ such that $|\lambda_{p_0}|=1.$ 

Since $F\in\mathcal{K},$ we know that $F$ satisfies the Maass relation:
\begin{equation}
A\left(F,\begin{pmatrix}m&r/2\\r/2&n\end{pmatrix}\right)=\sum_{d|(m,n,r)}d^{k-1}A\left(F,\begin{pmatrix}\frac{mn}{d^2}&\frac{r}{2d}\\ \frac{r}{2d}&1\end{pmatrix}\right)\label{5.1}
\end{equation}

Then the Maass relation gives:
\[A(F,p_0I_2)=A\left(F,\begin{pmatrix}p_0^{2}&\\&1\end{pmatrix}\right)+p_0^{k-1}A(F,I_2).\]
This gives:
\[\sum_{p\in\mathcal{P}}\lambda_pA(\mathbb{P}_{pI_2},p_0I_2)=\sum_{p\in\mathcal{P}}\lambda_pA\left(\mathbb{P}_{pI_2},\begin{pmatrix}p_0^2&\\&1\end{pmatrix}\right)+\sum_{p\in\mathcal{P}}\lambda_pp_0^{k-1}A(\mathbb{P}_{pI_2},I_2).\]
Therefore, we obtain
\begin{flalign*}
|\lambda_{p_0}A(\mathbb{P}_{p_0I_2},p_0I_2)|&\leq\sum_{\substack{p\in\mathcal{P}\\ p\neq p_0}}|\lambda_p||A(\mathbb{P}_{pI_2},p_0I_2)|\\
&+\sum_{p\in\mathcal{P}}|\lambda_p|\left|A\left(\mathbb{P}_{pI_2},\begin{pmatrix}p_0^2&\\&1\end{pmatrix}\right)\right|\\
&+\sum_{p\in\mathcal{P}}|\lambda_p|p_0^{k-1}|A(\mathbb{P}_{pI_2},I_2)|.
\end{flalign*}
Then combine Lemma \ref{3.1}, Lemma \ref{3.2} and Lemma \ref{3.3}, and we obtain:
\[\left\{1+o\left(\frac{1}{k^{100}}\right)\right\}|\lambda_{p_0}|=o\left(\frac{M^k}{k^{k\epsilon/2}}\right),\]
when $k$ is large. This shows that $|\lambda_{p_0}|<1$ when $k$ is large. A contradiction to the assumption that $|\lambda_{p_0}|=1.$
\qed

\section{Proof of Corollary \ref{2.2}}
\pf Here we need to decompose $\mathcal{B}_k$ into four disjoint sets, denoted by $\mathcal{B}_k^{i}$ and $i=1,2,3,4$. That is,
\[\mathcal{B}_k=\mathcal{B}_k^1\cup\mathcal{B}_k^2\cup\mathcal{B}_k^3\cup\mathcal{B}_k^4,\]
where
\begin{flalign*}
\mathcal{B}_k^1&=\{F\in\mathcal{B}_k|\mbox{$F$ is a Saito-Kurokawa lift and $A(F,I_2)\neq0.$}\}\\
\mathcal{B}_k^2&=\{F\in\mathcal{B}_k|\mbox{$F$ is not a Saito-Kurokawa lift and $A(F,I_2)\neq0.$}\}\\
\mathcal{B}_k^3&=\{F\in\mathcal{B}_k|\mbox{$F$ is a Saito-Kurokawa lift and $A(F,I_2)=0.$}\}\\
\mathcal{B}_k^4&=\{F\in\mathcal{B}_k|\mbox{$F$ is not a Saito-Kurokawa lift and $A(F,I_2)=0.$}\}\\
\end{flalign*}
Denote by $\mathcal{M}$ of $\mathcal{S}_k(\Gamma_2)$ spanned by $\mathcal{B}_k^3\cup\mathcal{B}_k^4.$ Let $F\in\mathcal{B}_k.$ By Andrianov's formula, 
\[\zeta(s+1/2)L(s+1/2,\chi_{-4})\sum_{m}\frac{A(F,mI_2)}{m^{s+k-3/2}}=L(s,\spin,F)A(F,I_2),\]
(see \cite{An1} Theorem 2.4.1 or \cite{An2} Theorem 4.3.16 with $l=a=1,$ and $\gamma=\eta=\chi=\triv$) we see that $F\in\mathcal{B}_k^3\cup\mathcal{B}_k^4$ implies that $F$ is orthogonal to the space $\mathcal{V}$ by \ref{1.1}. This will give:
\[\mathcal{V}\subseteq\mathcal{M}^{\perp}.\]
Here $\mathcal{M}^{\perp}$ is the orthogonal completment of $\mathcal{M}$ with respect to the Petersson inner product. Let $\mathcal{K}_0$ be the subspace of $\mathcal{K}$ spanned by $\mathcal{B}_k^1.$ The relation above can be refined to
\[\mathcal{V}+\mathcal{K}_0\subseteq\mathcal{M}^{\perp}.\]

\begin{rem}
Notice that this is only a sum of two spaces but not necessary a direct sum.
\end{rem}

This gives:
\begin{equation}
\#\mathcal{B}_k^1+\#\mathcal{B}_k^2=\#\{F\in\mathcal{B}_k|A(F,I_2)\neq 0\}\geq\dim(\mathcal{V}+\mathcal{K}_0)\geq\dim(\mathcal{W}+\mathcal{K}_0) .\label{2.5}
\end{equation}
Here $\mathcal{W}$ is the subspace of $\mathcal{V}$ which we found in Theorem \ref{2.1}. Since $\mathcal{W}\cap\mathcal{K}=\{0\},$ we know that $\mathcal{W}\cap\mathcal{K}_0=\{0\}.$ Then we obtain:
\[\#\mathcal{B}_k^1+\#\mathcal{B}_k^2\geq \dim(\mathcal{W}+\mathcal{K}_0)=\dim(\mathcal{W}\oplus\mathcal{K}_0)=\dim\mathcal{W}+\dim\mathcal{K}_0.\]
Since $\#\mathcal{B}_k^1=\dim\mathcal{K}_0$, we finally obtain:
\[\#\mathcal{B}_k^2\geq \dim\mathcal{W}.\]
\qed

\section{Loose Ends}
By virtual of Corollary \ref{2.2}, it suffices to find the lower bound for the dimension of $\mathcal{V},$ the subspace of $\mathcal{S}_k(\Gamma_2)$ spanned by $\{\mathbb{P}_{nI_2}|n\in\mathbb{N}\}.$ We can establish the following result:
\begin{theorem}
\upshape\KAI{Let $\mathcal{V}$ be the vector space spanned by $\mathbb{P}_{nI_2}.$ Let $\epsilon>0$. Then for any $\epsilon>0,$ we can find sufficiently large $k$ such that
\[\dim\mathcal{V}\geq \left(\frac{1}{2\sqrt{2}e\pi}-\epsilon\right)k.\]
}\label{6.1}
\end{theorem}
The proof could be quite similar to that of Theorem \ref{2.1}: for any $\epsilon>0,$ we set
\[\mathcal{Q}=\left\{n\in\mathbb{N} \left|n\leq  \left(\frac{1}{2\sqrt{2}e\pi}-\epsilon\right)k\right.\right\}.\]
Then we need to show that elements in the set $\{\mathbb{P}_{nI_2}|n\in\mathcal{Q}\}$ are linearly independent. Let $N$ be the largest element in $\mathcal{Q}.$

Similar to Theorem \ref{2.1}, we want to investigate the matrix
\[\mathcal{A}=\left(A(\mathbb{P}_{mI_2},nI_2)\right)_{m,n\in\mathcal{Q}}.\]
However, we need to normalize the Fourier coefficients of $\mathbb{P}_{nI_2}.$ That means, we define another matrix
\[\widetilde{\mathcal{A}}=\left(\tilde{a}_{m,n}\right)_{m,n\in\mathcal{Q}}\]
with
\[\widetilde{a}_{m,n}=\frac{\det(mI_2)^{l/2}}{\det(nI_2)^{l/2}}a_{m,n}=\frac{\det(mI_2)^{l/2}}{\det(nI_2)^{l/2}}A(\mathbb{P}_{mI_2},nI_2).\]
Notice that $\widetilde{\mathcal{A}}$ is invertible if and only if $\mathcal{A}$ is invertible since $\widetilde{\mathcal{A}}$ is equal to $\mathcal{A}$ multiplied by two (invertible) diagonal matrices on two sides.

By the idea of Lemma \ref{3.1}, we can prove:
\[\widetilde{a}_{m,n}=8\delta_{m,n}+O\left(\frac{4\sqrt{2}\pi cN}{l}\right)^l=8I_N+o\left(\frac{1}{k^{100}}\right).\]
when $k$ is sufficiently large. Here $c$ is the constant in (\ref{6.5}) and we chan choose $c$ quite close to $\frac{e}{2}$. The left part is similar to Theorem \ref{2.1}

To prove this, we need the following lemma, which can be considered as a refined version of Lemma \ref{3.4} in a special case:
\begin{lem}
\upshape\KAI{
Let $C$ be an invertible $2\times 2$ matrix. Then the smallest eigenvalue of $C^{-1}{}^tC^{-1}$ is $\leq 2||C||^{-2},$ where $||\cdot||$ denotes the matrix norm $||C||=\Tr(C{}^tC)$.
}\label{6.2}
\end{lem}
\pf This can be obtained by direct calculation.
\qed

However, this result is not quite interesting since Theorem \ref{6.1} fails to exclude Saito-Kurokawa lifts. Indeed, by the construction of Saito-Kurokawa lifts, Waldspurger's theorem and the Iwaniec-Sanark result \cite{IS}, (Or one can refer to Theorem 26.1 in \cite{IK}. The condition $k\equiv0 (\mo 4)$ is not essential in the theorem once the global root number is $1$.) we can show that  
\[\#\mathcal{B}_k^1=\#\{F\in\mathcal{B}_k|\mbox{$F$ is a Saito-Kurokawa lift and $A(F,I_2)\neq0.$}\}\geq\left(\frac{1}{12}-\epsilon\right)k.\]
This is greater than the lower bound in Theorem \ref{6.1}. So a natural question is, how to find better lower bounds for $\dim\mathcal{V}.$ 

\section*{Acknowledgement}
The author is grateful to the referee for a thorough reading of this paper and for helpful comments and suggestions, which lead to the improvements of the paper. In particular, the author would like to thank the referee for the valuable comments for Theorem \ref{2.1} and Corollary \ref{2.4}.

\end{document}